\documentclass[
]{ceurart}

\usepackage{amsfonts}

\usepackage{latexsym}
\usepackage{amssymb}
\usepackage{amsmath}
\usepackage{amsthm}
\usepackage[mathscr]{eucal}
\usepackage{graphicx}
\usepackage{hyperref}
\usepackage{caption}
\usepackage{subcaption}

\newtheorem{theorem}{Theorem}
\newtheorem{definition}{Definition}
\newtheorem{proposition}{Proposition}

\begin{document}

\copyrightyear{2021}
\copyrightclause{Copyright for this paper by its authors.
  Use permitted under Creative Commons License Attribution 4.0
  International (CC BY 4.0).}

\conference{Algorithms, Computing and Mathematics Conference (ACM 2021), Chennai, India}

\title{Two Applications of Graph Minor Reduction}

\author[1]{Ghurumuruhan Ganesan}[%
email=gganesan82@gmail.com,
]
\address[1]{Institute of Mathematical Sciences, HBNI, Chennai}

\begin{abstract}
In this paper, we study two applications of graph minor reduction. In the first part of the paper, we introduce a variant of the boxicity, called strong boxicity, where the rectangular representation satisfies an additional condition that each rectangle contains at least one point not present in any other rectangle. We show how the strong boxicity of a graph~\(G\) can be estimated in terms of the strong boxicity of a  minor~\(H\) and the number of edit operations needed to obtain~\(H\) from~\(G.\) In the second part of the paper, we consider false data injection (attack) in a flow graph~\(G\)
and quantify the subsequent effect on the state of edges of~\(G\) via the \emph{edge variation
factor}~\(\theta.\) We use minor reduction techniques to obtain bounds on~\(\theta\) in terms of the connectivity parameters of~\(G,\) when the attacker has complete knowledge of~\(G\) and also discuss stealthy attacks with partial knowledge of the flow graph.

\end{abstract}


\maketitle

\renewcommand{\theequation}{\thesection.\arabic{equation}}
\setcounter{equation}{0}
\section{Introduction} \label{intro}
Graph minor reduction is an important problem from both theoretical and application perspectives. In particular, graph minor reduction algorithms are extensively used in computer science to solve a variety of problems related to network routing and design. For a detailed survey of algorithmic aspects of graph minor reduction, we refer to the survey~\cite{bien}. 

In this paper, we study two theoretical applications of graph minor reduction in estimating the strong boxicity of an arbitrary graph and determining conditions for stealthy attacks on flow graphs. Below we discuss these two issues in that order.

\subsection*{Strong Boxicity}
The boxicity of a graph~\(G\)~\cite{rob} is the smallest integer~\(k\) such that~\(G\) admits a rectangular representation in~\(\mathbb{R}^k,\) where~\(\mathbb{R}\) denotes the real line. Boxicity as defined above is finite and is no more than the total number of vertices in~\(G.\)  Since then many bounds for boxicity has been obtained in terms of various graph parameters including treewidth~\cite{chan2} and maximum degree~\cite{chan, esp} and also in terms of poset dimension~\cite{adig, sco}. 

In the first part of the paper, we define a variant of boxicity which we call strong boxicity where we impose the additional condition that no rectangle corresponding to a particular vertex is covered by the rectangles corresponding to the other vertices. We find bounds for strong boxicity in terms of boxicity and estimate the strong boxicity of a general graph in terms of the strong boxicity of a minor and the number of edit operations needed to obtain the minor. 



\subsection*{Stealthy Attacks in Flow Graphs}
Flow graphs are expected to play an important role in the future as more and more systems are being automated through software. This also leads to potential vulnerabilities as software defined networks are themselves prone to attack and therefore it is important to study malicious data attacks on such automated flow networks. A typical example is that of a power grid where power flow through transmission lines is often subject to stealthy attacks (see for example,~\cite{kosut}).  The survey by~\cite{he} also describes various aspects of cyber-physical attacks and defence strategies on the smart grid, from a network layer perspective.

Flow graphs also frequently arise in the analysis of control systems~\cite{kuo,faw} where the node variables could be either electrical (like for example voltages, currents etc,) or mechanical (like position, angle etc.) in nature. Detection of false data injection (either unintentional or intentional) is crucial here as well and steps must be taken to ensure that the measurements are as authentic as possible.

In the second part of the paper, we are interested in studying how stealthy attacks affect the state of edges in a general flow graph. Different edges are affected differently and we quantify this effect via the \emph{edge variation factor}. In our main result Theorem~\ref{cap_thm} below, we estimate the maximum possible edge variation factor in terms of connectivity parameters of the flow graph, when the attacker has complete knowledge of the flow graph. In Section~\ref{prac_sec}, we also discuss possibility of attacks with partial flow graph knowledge.

The paper is organized as follows. In Section~\ref{sec_str_box}, we define the concept of strong boxicity and determine bounds for strong boxicity in terms of boxicity and also given examples of graphs with strong boxicity exactly equal to~\(2.\) Next in Section~\ref{sec_minor_red}, we show how strong boxicity of a graph can be estimated in terms of the strong boxicity of a minor and the number of edit operations needed to obtain the corresponding minor. In Section~\ref{sec_stealth}, we state and prove our main result  Theorem~\ref{cap_thm} regarding existence of stealthy attacks in general flow graphs and in Section~\ref{prac_sec}, we discuss stealthy attacks in the presence of partial information.


\renewcommand{\theequation}{\thesection.\arabic{equation}}
\setcounter{equation}{0}
\section{Strong boxicity}\label{sec_str_box}
Let~\(\mathbb{R}\) denote the real line and for integer~\(d \geq 1\) define a rectangle~\(A \subset \mathbb{R}^d\) to be a closed set of the form~\(\prod_{i=1}^{d} [a_i,b_i] \subset \mathbb{R}^{d},\) where~\(a_i \leq b_i\) are real numbers. We define~\(A^{0} := \prod_{i=1}^{n} (a_i, b_i)\) to be the interior of the rectangle~\(A\) where~\((a_i,b_i)\) denotes the open interval with endpoints~\(a_i\) and~\(b_i\) and Set~\(\partial A := A \setminus A^{0}\) to be the boundary of the rectangle~\(A.\) Also for~\(y > 0, x \in \mathbb{R}^d,\) we let~\(B_d(x,y) := x + \left[-\frac{y}{2},\frac{y}{2}\right]^d\) be the~\(d-\)dimensional square of side length~\(y\) centred at~\(x.\)

\begin{figure}
\centering
\begin{subfigure}{0.5\textwidth}
\centering
\includegraphics[width=2in, trim= 50 180 50 40, clip=true]{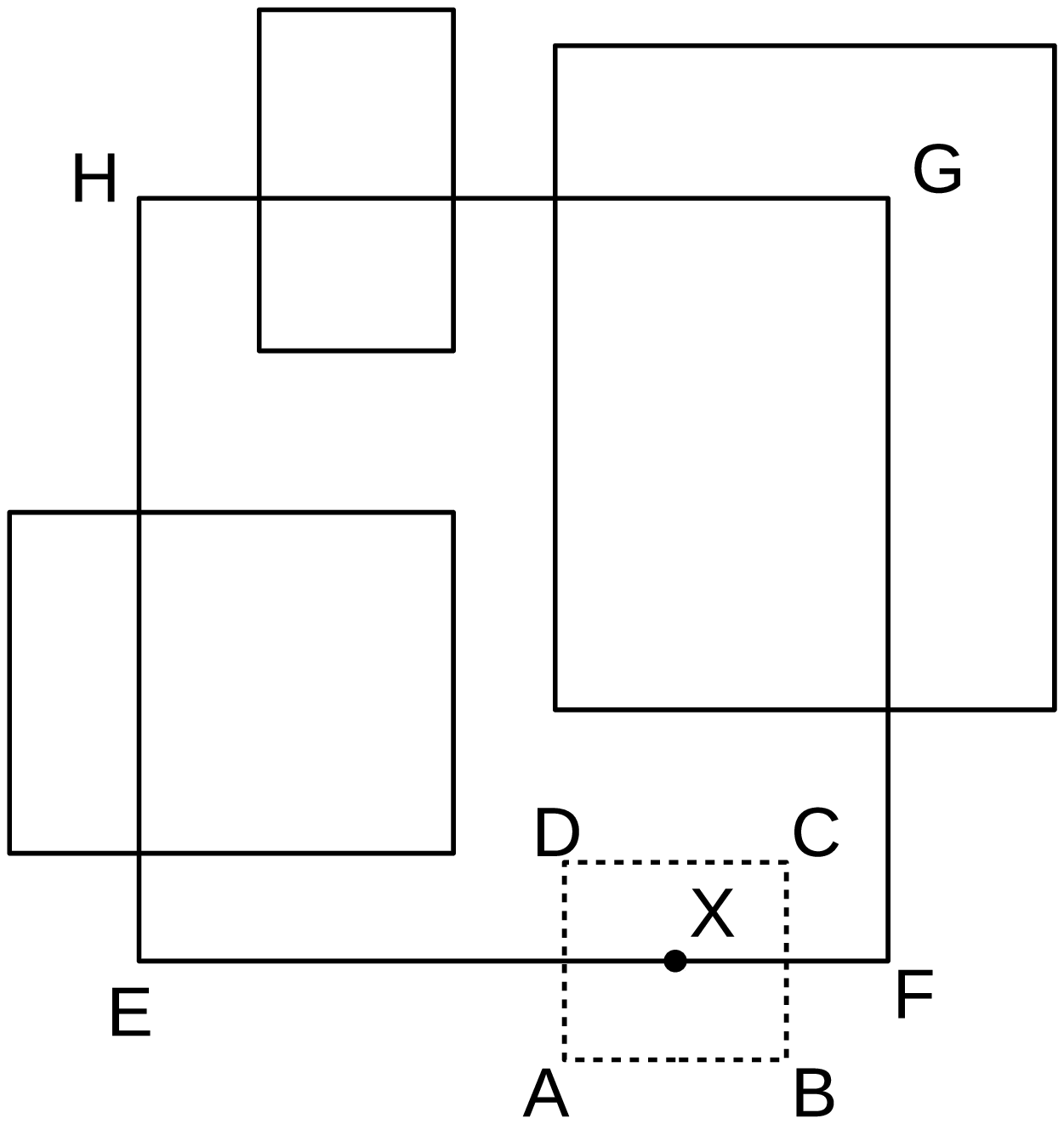}
  \caption{} 
\end{subfigure}
\begin{subfigure}{0.5\textwidth}
\centering
   \includegraphics[width=2in, trim= 50 180 50 40, clip=true]{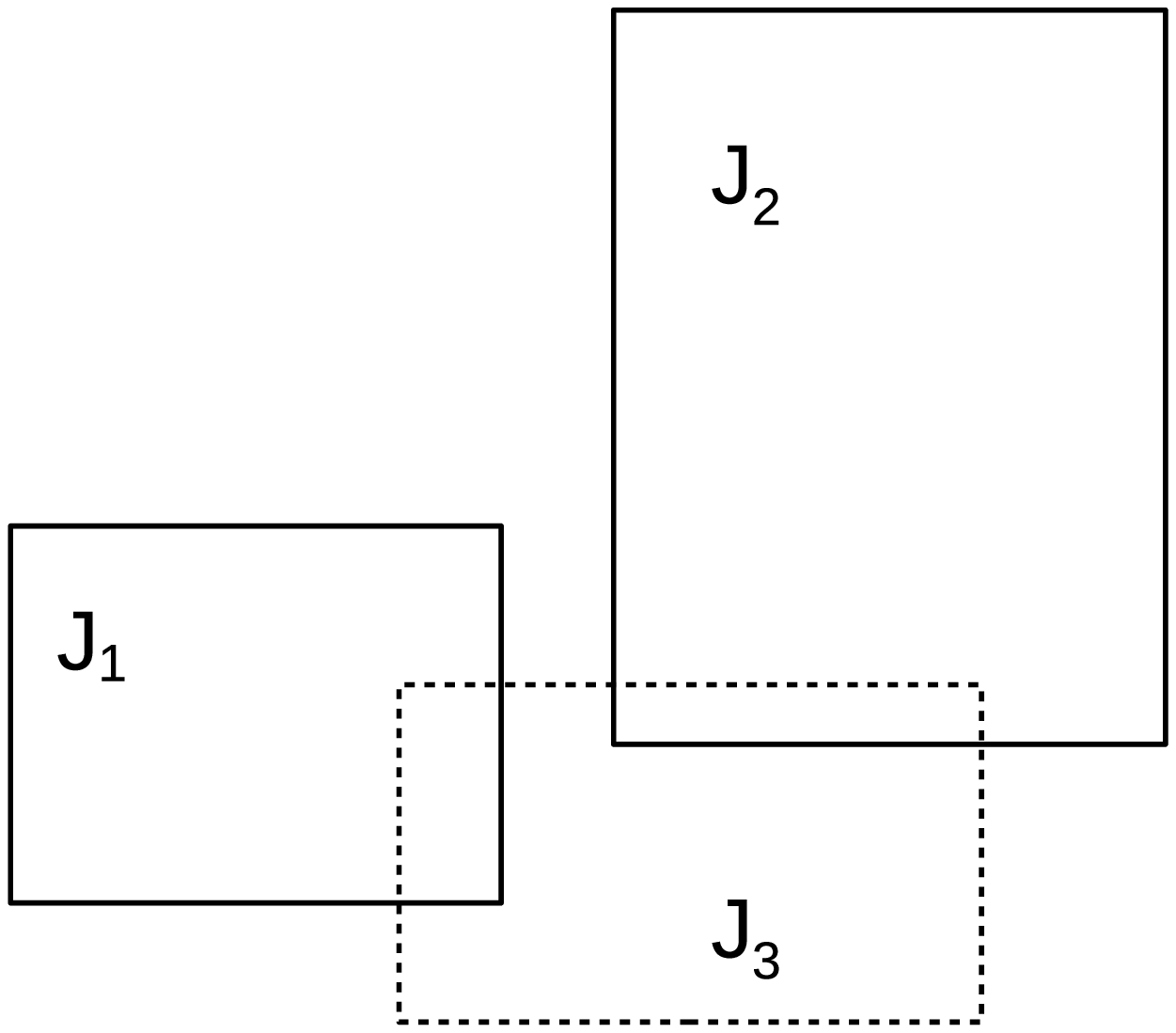}
   \caption{} 
  \end{subfigure}

\caption{\((a)\) The condition~\((c2)\) ensures that each rectangle has an exclusive boundary point and a small surrounding neighbourhood. \((b)\) Strong rectangular representation of the graph~\(H\) with vertex set~\(\{1,2,3\}\) and edge set~\(\{(1,2),(2,3)\}.\)}
\label{diag_squares}
\end{figure}

Let~\(G = (V,E)\) be a graph with vertex set~\(V = \{1,2,\ldots,n\}\) and define~\((i,j) \in E\) to be the edge with endvertices~\(i\) and~\(j.\)
\begin{definition}\label{str_box_def}
We say that a set of rectangles~\(\{I_i\}_{1 \leq i \leq n}\) in~\(\mathbb{R}^d, d \geq 1\) is a \emph{strong rectangular representation} of~\(G\) if the following two conditions hold:\\
\((c1)\) For any~\(1 \leq i \neq j \leq n,\) the rectangles
\begin{equation}\label{int_cond}
I_i \cap I_j \neq \emptyset \text{ if and only if } (i,j) \in E.
\end{equation}
\((c2)\) For every~\(1 \leq v \leq n,\) there exists a point~\(x_v \in \partial I_v\) and a number~\(a_v >0\) such that
\begin{equation}\label{str_rep}
B_d(x_v,a_v) \subset \left(\bigcup_{1 \leq u \neq v \leq n} I_u\right)^c.
\end{equation}
\end{definition}
In other words, we prefer that no rectangle has its boundary covered completely by the remaining rectangles. The \emph{strong boxicity} of~\(G\) denoted by~\(s(G)\) is defined to be the smallest integer~\(d\) such that both the conditions~\((c1)-(c2)\) hold. In Figure~\ref{diag_squares}, we illustrate condition~\((c2)\) with a example involving a strong rectangular representation in~\(\mathbb{R}^2.\) The solid rectangle~\(ABCD\) is~\(I_v,\) the point~\(X = x_v\) and the dotted rectangle corresponds to~\(B_d(x_v,a_v).\)

We recall that the boxicity~\(G\) denoted by~\(b(G)\) is defined to be smallest integer~\(d\) such that condition~\((c1)\) alone holds~\cite{rob}. By definition we therefore have~\(b(G) \leq s(G).\) To see that strong boxicity is not always equal to boxicity, we consider the graph~\(H\) with vertex set~\(\{1,2,3\}\) and edge set~\(\{(1,2),(2,3)\}.\) Setting~\(J_1 = [0,2], J_2 = [1,4]\) and~\(J_3 = [2,5],\) we see that condition~\((c1)\) in definition~\ref{str_box_def} holds and so~\(b(H) = 1.\) However, condition~\((c2)\) does not hold and so~\(s(H) \geq 2.\) Considering~\(J_1,J_2\) and~\(J_3\) to be squares such that~\(J_1\) and~\(J_3\) are disjoint and~\(J_2\) intersects both~\(J_1\) and~\(J_3\) partially, as shown in Figure~\ref{diag_squares}\((b),\) we get that~\(s(H)=2.\)

From the discussion in the above paragraph, we deduce that any connected graph containing at least two edges must have strong boxicity at least two. In the following result, we give examples of graphs whose strong boxicity is exactly~\(2\) and also obtain bounds for the strong boxicity in terms of the boxicity. We begin with some definitions. A tree is a connected graph containing~\(n\) vertices and~\(n-1\) edges for some~\(n \geq 2.\) A \emph{clique} in a graph~\(G\) is a complete subgraph of~\(G.\) We say that~\({\cal I}\) is a \emph{stable} set in~\(G\) if no two vertices are adjacent in~\(G.\) The graph~\(G\) with vertex set~\(\{1,2,\ldots,n+r\}\) is called a \emph{threshold graph}~\cite{maha} if~\(G\) contains a clique~\(K\) with vertex set~\(\{1,2,\ldots,n\}\) and a stable set~\(I = \{n+1,\ldots,n+r\}\) satisfying the following nested neighbourhood property: If~\({\cal N}(v)\) is the set of all neighbours of~\(v\) in the graph~\(G,\) then
\begin{equation}\label{nest_neigh}
{\cal N}(n+1) \supseteq {\cal N}(n+2) \supseteq \ldots \supseteq {\cal N}(n+r).
\end{equation}
\begin{proposition}\label{prop1} If~\(G\) is a tree or a threshold graph containing at least three vertices, then the strong boxicity~\(s(G) = 2.\)
In general, for any graph~\(G\) we have that
\begin{equation}\label{sbb}
b(G) \leq s(G) \leq b(G)+2.
\end{equation}
\end{proposition}
From~(\ref{sbb}), we see that any bound for boxicity can also be used to estimate the strong boxicity. Therefore using~\(b(G) \leq \min(m,n)\) where~\(m\) is the number of edges of~\(G\)~\cite{rob}, we get from~(\ref{sbb}) that~\(s(G) \leq \min(m,n)+2.\)


\emph{Proof of Proposition~\ref{prop1}}: First, we see that any tree or threshold graph containing at least three vertices has strong boxicity of exactly two.\\
\underline{\emph{Trees}}: We prove by induction on~\(n,\) the number of vertices in the tree. For~\(n=3\) the tree~\(T_3\) contains  two edges, say~\((1,2)\) and~\((2,3).\) To see that~\(s(T_3) = 2,\) we define~\(J_1,J_2\) and~\(J_3\) to be squares such that~\(J_1\) and~\(J_3\) are disjoint and~\(J_2\) intersects both~\(J_1\) and~\(J_3\) partially, as shown in Figure~\ref{diag_squares}\((b)\).

Now, let~\(T_{n+1}\) be a tree with vertex set~\(\{1,2,\ldots,n+1\}\) and suppose the vertex~\(n+1\) is a leaf attached to the vertex~\(v \in T_n := T_{n+1} \setminus \{n+1\}.\) The tree~\(T_n\) has strong boxicity two by induction assumption and so there exists a strong rectangular representation~\(\{J_i\}_{1 \leq i \leq n} \subset \mathbb{R}^2\) of~\(T_n.\) Moreover, there exists a vertex~\(x_v \in J_v\) and a real number~\(a_v> 0\) satisfying~(\ref{str_rep}). Setting~\(J_{n+1} := B_{2}\left(x_v,\frac{a_v}{2}\right),\) we get that~\(\{J_i\}_{1 \leq i \leq n+1}\) forms a strong rectangular representation of~\(T_{n+1}.\)

This is illustrated in Figure~\ref{diag_squares}\((a),\) where~\(J_v = EFGH\) intersects other rectangles in the rectangular representation. However, the point~\(X = x_v\) and a small surrounding neighbourhood is unique to~\(J_v.\) The rectangle~\(J_{n+1} = ABCD\) intersecting only~\(J_v\) is shown in dotted lines.

\underline{\emph{Threshold graphs}}: Let~\(G = (K,I)\) be any threshold graph with\\\(K = \{1,2,\ldots,n\}\) being a clique of size~\(n\) and~\(I = \{n+1,\ldots,n+r\}\) being a stable set. Because of the nested neighbourhood property~(\ref{nest_neigh}), we assume that~\(N(n+i) = \{1,2,\ldots,l_i\}\) for some~\(n \geq l_1 \geq l_2 \geq \ldots  \geq l_r \geq 1.\)

\begin{figure}[tbp]
\centering
\includegraphics[width=2in, trim= 50 180 50 10, clip=true]{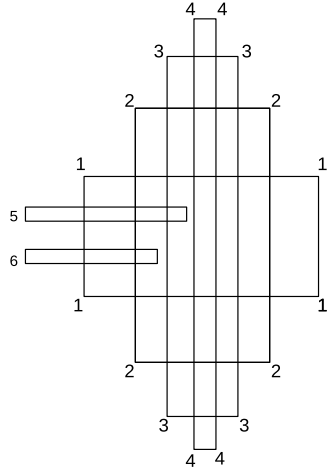}
\caption{Strong rectangular representation of a threshold graph.}
\label{thresh_box_fig}
\end{figure}

For~\(1 \leq i \leq n,\) let~\(J_i\) be the~\(\frac{1}{i} \times i\) rectangle in~\(\mathbb{R}^2\) centred at the origin. The rectangles~\(\{J_i\}_{1 \leq i \leq n}\) form a strong representation of~\(K.\) We then let~\(J_{n+1},\ldots,J_{n+r}\) be small disjoint rectangles with the property that~\(J_{n+i}\) intersects only the rectangles~\(J_1,\ldots,J_{l_i}\) and no other rectangle; this is illustrated in Figure~\ref{thresh_box_fig}  for the case~\(K = \{1,2,3,4\}\) and~\(I = \{5,6\}\) having neighbourhoods~\({\cal N}(5) = \{1,2,3\}\) and~\({\cal N}(6) = \{1,2\}.\) The rectangles with corners labelled~\(i, 1 \leq i \leq 4\) form the strong rectangular representation of~\(K.\) The rectangles labelled~\(5\) and~\(6\) are~\(J_5\) and~\(J_6,\) respectively. This completes the proof that any threshold graph containing at least three vertices has a strong boxicity of exactly two.

In the rest of the proof we prove~(\ref{sbb}). We use the following boundary relation throughout. For integer~\(d \geq 1\) let~\(U_d = \prod_{i=1}^{d}[a_i,b_i] \subset \mathbb{R}^d\) be any rectangle and let~\(\partial U_d\) be its boundary.
We then have that
\begin{equation}\label{par_un}
\partial U_d \supseteq \partial U_{d-1} \times [a_d,b_d].
\end{equation}
Indeed, writing~\(\partial U_d = U_d \setminus U^{0}_d = U_{d-1} \times [a_d,b_d] \setminus U^{0}_{d-1} \times (a_d,b_d)\)
and using
\begin{eqnarray}
U_{d-1} \times [a_d,b_d] &=& \left(U_{d-1}^{0} \cup \partial U_{d-1} \right)\times [a_d,b_d] \nonumber\\
&\supseteq& U_{d-1}^{0} \times (a_d,b_d) \bigcup \partial U_{d-1} \times [a_d,b_d], \nonumber
\end{eqnarray}
we get~(\ref{par_un}).

To prove~(\ref{sbb}), we let~\(\{I_i\}_{1 \leq i \leq n} \subset \mathbb{R}^{b}, b = b(G)\) be a rectangular representation of~\(G.\) For~\(1 \leq i \leq n,\) we define the rectangle~\(J_i := I_i \times [0,i] \times [i,n]\) and show first that~\(\{J_i\}_{1 \leq i \leq n} \subset \mathbb{R}^{b+2}\) forms a rectangular representation of~\(G.\) Indeed if~\((i,j)\) is an edge in~\(G\) then~\(I_i \cap I_j \neq  \emptyset\) and so
\begin{equation}\nonumber
J_i \cap J_j = (I_i \cap I_j) \times [0,\min(i,j)] \times [\max(i,j),n] \neq \emptyset.
\end{equation}

To see that the rectangles~\(\{J_i\}\) form a strong rectangular representation of the graph~\(G,\) we let~\(v \in G\) be any vertex and let~\(x_v \in \partial I_v\) be any point in the boundary of~\(I_v.\) Setting~\(y_v:=  (x_v,v,v)\) we have  from~(\ref{par_un}) that~\(y_v \in \partial J_v.\) Also~\(y_v \notin J_u\) for any~\(u \neq v.\) Letting~\(a_v := \frac{1}{4},\) we also get that no point of~\(B_{b+2}(y_v,a_v)\) belongs to~\(\bigcup_{1 \leq i \neq v \leq n }J_i.\) Therefore~\(s(G)\leq b(G)+2\) and this proves~(\ref{sbb}).~\(\qed\)






\setcounter{equation}{0}
\renewcommand\theequation{\thesection.\arabic{equation}}
\section{Minor Reduction}\label{sec_minor_red}
In this subsection, we estimate the strong boxicity of a graph~\(G\) by ``converting"~\(G\) into a graph~\(H\) with small strong boxicity through appropriate edit operations. Formally, a \emph{deletion operation} on~\(G\) is either a vertex deletion or edge deletion and a contraction operation is defined as follows. Let~\(e = (u,n)\) be an edge in~\(G\) for some~\(1 \leq u \leq n-1.\) The graph~\(H\) obtained by \emph{contracting} the edge~\(e\) has vertex set~\(\{1,2,\ldots,n-1\}\) and still has all edges not containing~\(n\) as an endvertex. In addition, in the graph~\(H,\) the vertex~\(u\) is now adjacent to every vertex of~\(N_G(n) \setminus \{u\}.\) Henceforth, we denote a deletion operation or a contraction as an \emph{edit operation}.

The graph~\(H\) obtained from an edit operation on~\(G\) as described above, is called a \emph{minor} of~\(G.\) We say that a minor~\(H\) is obtained from~\(G\) after~\(a_1\) vertex deletions,~\(a_2\) edge deletions and~\(a_3\) contractions if there are graphs\\\(G = H_0,H_1,H_2,\ldots,H_r = H, r = a_1+a_2+a_3\) such that~\(H_i\) is obtained by either a vertex deletion, edge deletion or contraction of~\(H_{i-1}\) and the total number of vertex deletions is~\(a_1\) and the total number of edge deletions is~\(a_2.\) We have the following result regarding the strong boxicity.
\begin{theorem}\label{box_thm} If~\(H\) is a minor that is obtained from~\(G\) after~\(\alpha_v\) vertex deletions and~\(\alpha_e\) edge deletions,
then
\begin{equation}\label{box_eq}
s(H) - \alpha_e \leq s(G) \leq s(H) + \alpha_v + \alpha_e.
\end{equation}

If~\(F\) is a minor that is obtained from~\(G\) after~\(\beta_v\) vertex deletions,~\(\beta_e\) edge deletions and~\(\beta_c\) contractions,
then
\begin{equation}\label{box_eq2}
s(G) \leq s(F) + \beta_v + \beta_e + 2 \beta_c.
\end{equation}
Moreover both~(\ref{box_eq}) and~(\ref{box_eq2}) hold with strong boxicity~\(s(.)\) replaced by boxicity~\(b(.).\)
\end{theorem}
In practice, we choose~\(H\) to be a known graph with low strong boxicity. As an example, suppose~\(G\) is a connected graph on~\(n\) vertices having~\(n+r\) edges. We pick and remove~\(r+1\) edges from~\(G\) to get a tree~\(H,\) whose strong boxicity is two by proposition~\ref{prop1}. From relation~(\ref{box_eq}) in theorem~\ref{box_thm}, we therefore get that~\(s(G) \leq r+3.\)

It suffices to prove Theorem~\ref{box_thm} for a single edit operation and we consider the three cases vertex deletion, edge deletion and edge contraction separately below. Also we prove for strong boxicity throughout and an analogous analysis holds for boxicity.

\underline{\emph{Vertex deletion}}: We show that for any vertex~\(v \in G,\) the strong boxicity
\begin{equation}\label{case_a}
s\left(F_v\right) \leq s(G) \leq s\left(F_v\right) +1,
\end{equation}
where~\(F_v = G \setminus \{v\}\) is the graph obtained by removing the vertex~\(v \in G.\) We prove for~\(v = n\) and an analogous analysis holds for every other vertex. If~\({\cal I} =\{I_i\}_{1 \leq i \leq n} \subset \mathbb{R}^{s}, s = s(G)\) is a strong rectangular representation of~\(G\) satisfying~(\ref{int_cond}), then~\({\cal I}\setminus \{I_n\}\) is a strong rectangular representation of~\(G \setminus \{n\}\) and so the first inequality in~(\ref{case_a}) is true.

For the second inequality, we let~\(N(n)\) be the neighbours of~\(n\) in~\(G\) and let~\(\{J_i\}_{1 \leq i \leq n-1} \subset \mathbb{R}^{k}, k = s\left(F_n\right)\) be a strong rectangular representation of~\(F_n\) with~\(J_i = \prod_{l=1}^{k}[a_{i,l},b_{i,l}].\) Setting~\(M_{up} := \max_{i,l}b_{i,l}\) and~\(M_{low} = \min_{i,l} a_{i,l}\) we define the rectangles~\(\{I_w\}_{1 \leq w \leq n} \subset \mathbb{R}^{k+1}\) as
\begin{equation}\label{iw_def}
I_w :=
\left\{
\begin{array}{cc}
J_w \times [0,3] & \text{ for } w \notin \{n\} \cup N(n) \\
J_w \times [2,5] & \text{ for } w \in N(n) \\
\left[M_{low}, 3M_{up}\right]^{k} \times [4,6] & \text{ for } w = n.
\end{array}
\right.
\end{equation}
By construction, the rectangles~\(\{I_i\}_{1 \leq i \leq n}\) form a rectangular representation of~\(G.\) We now use~(\ref{par_un}) to see that~\(\{I_i\}\) also form a strong rectangular representation of~\(G.\) Indeed for~\(w \notin \{n\} \cup N(n),\) we let~\(x_w \in \partial J_w\) and~\(a_w> 0\) be such that
\begin{equation}\label{str_rep2}
B_k(x_w,a_w) \subset \left(\bigcup_{1 \leq u \neq w \leq n} J_u\right)^c.
\end{equation}
From the first line in~(\ref{iw_def}) and~(\ref{par_un}) we get that~\(y_w := (x_w,0) \in \partial I_w\) and from~(\ref{str_rep2}) and the second and third lines in~(\ref{iw_def}) we further get
\begin{equation}\label{str_rep3}
y_w + B_k(0,a_w) \times [-1,1] \subset \left(\bigcup_{1 \leq u \neq w \leq n} I_u\right)^c.
\end{equation}
Therefore choosing~\(a_w\) smaller if necessary we get
\begin{equation}\label{str_rep4}
B_{k+1}(y_w,a_w) \subset \left(\bigcup_{1 \leq u \neq w \leq n} I_u\right)^c.
\end{equation}

If~\(w \in N(n)\) then letting~\(x_w\) be as in~(\ref{str_rep2}), we choose~\(y_w = (x_w,2).\) Arguing as before and choosing~\(a_w\) smaller if necessary, we get~(\ref{str_rep4}). Finally, if~\(w = n,\) then we choose~\(y_w\) to be the~\((k+1)\)-tuple~\((3M_{up},\ldots,3M_{up},6)\) so that~\(y_w \in \partial I_w.\) Setting~\(a_w = M_{up}\) and using the definition of~\(M_{up}\) we again get~(\ref{str_rep4}). Thus the strong boxicity~\(s(G) \leq k+1.\)

\underline{\emph{Edge deletion}}: For any edge~\(e = (u,v) \in G\) we show that
\begin{equation}\label{case_b}
s\left(H_e\right) -1 \leq s(G) \leq s\left(H_e\right) +1,
\end{equation}
where~\(H_e = G \setminus \{e\}\) is the graph obtained by removing the edge~\(e.\)


To prove the lower bound in~(\ref{case_b}), we let~\(\{I_i\}_{1 \leq i \leq n} \subset \mathbb{R}^{s}, s = s(G)\) be a strong rectangular representation of~\(G.\)
We define the rectangles~\(\{A_i\}_{1 \leq i \leq n} \subset \mathbb{R}^{s+1}\) as follows:
\begin{equation}\nonumber
A_i :=
\left\{
\begin{array}{cc}
I_i \times [0,4] & \text{ for } i \neq u,v \\
I_i \times [1,2] & \text{ for } i=u \\
I_i \times [3,5] & \text{ for } i = v.
\end{array}
\right.
\end{equation}
Arguing as in the vertex deletion case, we get that the rectangles~\(\{A_i\}_{1 \leq i \leq n}\) form a strong rectangular representation of~\(H_e\) and so the strong boxicity~\(s\left(H_e\right) \leq s+1.\)

For the upper bound in~(\ref{case_b}), we use the fact that the graph~\(H_e\) has~\(n\) vertices and let~\(L_1,\ldots,L_n \subset \mathbb{R}^{r}, r = s(H_e)\) be a strong rectangular representation of~\(H_e.\) As before if~\(L_i = \prod_{l=1}^{r} [a_{i,l},b_{i,l}]\)
then we set~\(W_{up} := \max_{i,l} b_{i,l}\) and~\(W_{low} := \min_{i,l}a_{i,l}.\)

Letting~\(N_e(v)\) be the neighbours of~\(v\) in~\(H_e,\) we now define the rectangles~\(\{R_i\}_{1 \leq i \leq n} \subset \mathbb{R}^{r+1}\) as follows:
\begin{equation}\nonumber
R_i :=
\left\{
\begin{array}{cc}
L_i \times [0,3] & \text{ for } i \notin \{u,v\} \cup N_e(v) \\
L_i \times [2,5] & \text{ for } i \in \{u\} \cup N_e(v) \\
\left[W_{low}, 3W_{up}\right]^{r} \times [4, 6] & \text{ for } i = v.
\end{array}
\right.
\end{equation}
Arguing as in the vertex deletion case, we get that the rectangles~\(\{R_i\}_{1 \leq i \leq n}\) form a strong rectangular representation of~\(G\) and so the strong boxicity~\(s(G) \leq r+1.\)

\emph{\underline{Edge contraction}}: We now consider the remaining case where~\(G_e\) is the graph obtained by the contracting
the edge~\(e =(u,n) \in G\) to the vertex~\(u \in G_e\) and show that
\begin{equation}\label{case_c}
s(G) \leq s(G_e)+2.
\end{equation}

The graph~\(G_e\) has~\(n-1\) vertices and we
let~\(\{S_i\}_{1 \leq i \leq n-1} \subset \mathbb{R}^{q}, q = s(G_e),\) be a strong rectangular representation of~\(G_e.\) If~\(N(u)\) and~\(N(n)\) are the neighbours of~\(u\) and~\(n\) respectively in~\(G,\) then~\(n  \in N(u)\) and~\(u \in N(n)\) and we define the rectangles~\(\{T_i\}_{1 \leq i \leq n} \subset \mathbb{R}^{q+2}\) as follows:
\begin{equation}\label{t_def}
T_i :=
\left\{
\begin{array}{cc}
S_u \times [0,6] \times [3,7] & \text{ for } i = u\\
S_u \times [4,10] \times [6,10] & \text{ for } i = n\\
S_i \times [0,10] \times [0,5] & \text{ for } i \in N(u) \setminus \left(\{n\} \cup N(n)\right) \\
S_i \times [8,10] \times [0,10] & \text{ for } i \in N(n) \setminus \left(\{u\} \cup N(u)\right)\\
S_i \times [0,10] \times [0,10] & \text{ otherwise }.
\end{array}
\right.
\end{equation}

We first argue that the rectangles~\(\{T_i\}_{1 \leq i \leq n}\) form a rectangular representation of~\(G.\)
Let~\((i,j)\) be an edge not in~\(G.\) If~\((i,j)\) is not of the form~\(i =u ,j \in N(n)\setminus N(u)\)
or~\(i=n,j \in N(u) \setminus N(n),\) then the edge~\((i,j)\) is not present in~\(G_e\) as well and so~\(S_i \cap S_j = \emptyset.\)
Consequently~\(T_i \cap T_j = \emptyset.\)

If~\( i=u\) and~\(j \in N(n) \setminus N(u),\)
then by definition of contraction
we have that~\((u,j) \in G_e\) and so~\(S_u \cap S_j \neq \emptyset.\) But from
the first and fourth lines in~(\ref{t_def}), we get that~\(T_u \cap T_j = \emptyset.\)
Similarly if~\(i=n\) and~\(j \in N(u) \setminus N(n),\)
then~\(S_i \cap S_j = S_u \cap S_j \neq \emptyset.\) But from the second and third lines in~(\ref{t_def})
we get that~\(T_i \cap T_j = \emptyset.\)

Next let~\((i,j)\) be an edge present in~\(G.\) If~\(i,j \neq u\) or~\(n,\) then~\((i,j) \in G_e\) and so~\(S_i \cap S_j \neq \emptyset.\)
From the last three lines of~(\ref{t_def}) we then get that~\[T_i \cap T_j \supseteq (S_i \cap S_j) \times [8,10] \times [0,5] \neq \emptyset.\]
If~\(i=u\) and~\(j=n \in N(u),\) then from the first two lines of~(\ref{t_def}) we get that
\begin{equation}\label{tun}
T_u \cap T_n = S_u \times [4,6] \times [6,7] \neq \emptyset.
\end{equation}
If~\(i=u\) and~\(j \in N(u)\setminus (\{n\} \cup N(n))\) then~\(S_u \cap S_j \neq \emptyset\) and so from the first and third lines in~(\ref{t_def}),
we get that~\[T_u \cap T_j = (S_u \cap S_j) \times [0,6] \times [3,5] \neq \emptyset.\]

Finally if~\(i=n\) and~\(j \in N(n)\setminus \{u\} \) then using the fact that~\(S_u \cap S_j \neq \emptyset\) and the second and fourth lines in~(\ref{t_def}), we get that~\[T_n \cap T_j \supset (S_u \cap S_j) \times [8,10] \times [6,10] \neq \emptyset.\] This proves the rectangles~\(\{T_i\}\) for a rectangular representation of~\(G.\)

To see that~\(\{T_i\}\) in fact form a strong rectangular representation, it suffices to see that~(\ref{str_rep}) in condition~\((c2)\) holds for~\(v = u\) and~\(v = n.\) Using the fact that~\(\{S_i\} \subset \mathbb{R}^s\) form a strong rectangular representation of~\(G_e,\) we let~\(x_u \in \mathbb{R}^{s}\) and~\(a_u > 0\) be such that
\begin{equation}\nonumber
B_q(x_u,a_u) \subset \left(\bigcup_{1 \leq i \neq u \leq n-1} S_i\right)^c.
\end{equation}
Using the first and second lines of~(\ref{t_def}), we then set~\(y_u := (x_u,0,3)\) and~\(y_n:= (x_u,10,6)\) respectively and choose~\(a_u\) smaller if necessary to get
\begin{equation}\nonumber
B_{q+2}(y_u,a_u) \subset \left(\bigcup_{1 \leq i \neq u \leq n} T_i\right)^c \text{ and } B_{q+2}(y_n,a_u) \subset \left(\bigcup_{1 \leq i \leq n-1} T_i\right)^c.
\end{equation}
By~(\ref{par_un}), we have that~\(y_u \in \partial T_u\) and~\(y_n \in \partial T_n\) and so~\(s(G_e) \leq q +2.\)

\begin{figure}[tbp]
\centering
\includegraphics[width=2.5in, trim= 10 220 15 100, clip=true]{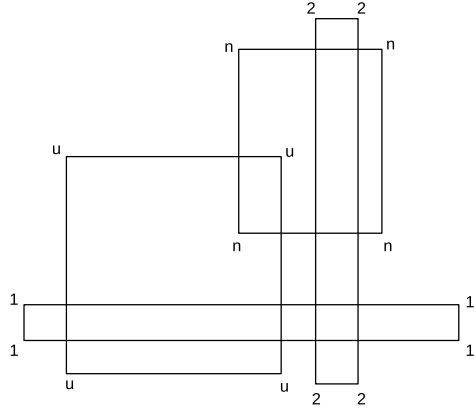}
\caption{Rectangles used in constructing a strong rectangular representation for~\(G.\)}
\label{cont_fig}
\end{figure}

Finally, in Figure~\ref{cont_fig}, we pictorially represent a possible choice for the rectangles in~\(\mathbb{R}^2\) that can be used to construct a strong rectangular representation for~\(G\) as follows: The rectangle~\(L_x\) with corners labelled~\(x \in \{u,n\}\) is used for the vertex~\(x\) and so we set~\(T_x = S_x \times L_x.\) Similarly, the rectangle~\(L_1\) with corners labelled~\(1\) is for neighbours of~\(u\) that are not adjacent to~\(n\) and the rectangle~\(L_2\) with corners labelled~\(2\) is for neighbours of~\(n\) that are not adjacent to~\(v.\) For the rest of the vertices, we pick any rectangle~\(L_0\) that is adjacent all the rectangles in Figure~\ref{cont_fig}. Defining the appropriate rectangles~\(T_i, 1 \leq i \leq n,\) we then get a strong rectangular representation of~\(G\) in~\(\mathbb{R}^{b+2}.\)~\(\qed\)

\emph{Proof of Theorem~\ref{box_thm}}: Follows from~(\ref{case_a}),~(\ref{case_b}) and~(\ref{case_c}).~\(\qed\)




%





\renewcommand{\theequation}{\thesection.\arabic{equation}}
\setcounter{equation}{0}
\section{Stealthy Attacks on Flow Graphs}\label{sec_stealth}
A flow graph is a graph~\(G = (V,E)\) with vertex set~\(V = \{1,2,\ldots,n\}\) and edge set~\(E\) (which we index as~\(\{n+1,\ldots,t\}\)) along with the following additional properties: The state of the system is given by a real valued vector~\(\mathbf{x} = [x_1,\ldots,x_n]^T\) where~\(x_i\) denotes the state of vector~\(i.\) An edge with index~\(f\) joining vertices~\(i\) and~\(j\) is assigned a \emph{gain}~\(B_{i,j} = B_{j,i} > 0\) and the \emph{flow} through the edge~\(f\) from~\(j\) to~\(i\) is given by
\begin{equation}\label{line_z}
B_{i,j}(x_i-x_j) = \mathbf{h}_f\cdot \mathbf{x},
\end{equation}
where~\(\mathbf{h}_f\) is the~\(1 \times n\) vector with exactly two non-zero entries:~\(B_{i,j}\) in~\(i^{th}\) position and~\(-B_{i,j}\) in the~\(j^{th}\) position. The \emph{net} flow into the vertex~\(i\) equals the sum of flows from all edges connected to~\(i\) and is given by
\begin{equation}\label{bus_z}
\sum_{j \sim i} B_{i,j}(x_i-x_j) = x_i \sum_{j \sim i} B_{i,j} - \sum_{j \sim i} B_{i,j} x_j = \mathbf{h}_i \cdot \mathbf{x},
\end{equation}
with~\(j \sim i\) denoting that vertices~\(j\) and~\(i\) are connected by an edge in~\(G.\)  The~\(1 \times n\) vector~\(\mathbf{h}_i\) has values~\(- B_{i,j}\) for positions~\(1 \leq j \neq i \leq n, j \sim i\) and has the value~\(H_{i,i} = \sum_{u \sim i}B_{i,u}.\) All other entries of~\(\mathbf{h}_i\)
are zero. Letting~\(\mathbf{H} = [\mathbf{h}_1^T,\ldots,\mathbf{h}_t^T]^T\) be the~\(t \times n\) gain matrix, we then have
\begin{equation}\label{h_sat}
\sum_{j=1}^{n} H_{i,j} = 0
\end{equation}
for all~\(1 \leq i \leq t.\)

Given the flow vector~\(\mathbf{H} \cdot \mathbf{x},\) the gain matrix~\(\mathbf{H}\) and a reference state (say~\(x_1\)), it is possible to calculate the overall state~\(\mathbf{x}:\) We first use~(\ref{line_z}) to  get~\(x_i-x_j\) across every edge~\((i,j) \in E.\) We then calculate the states of all vertices adjacent to the vertex~\(1\) and then iteratively calculate the states of the rest of the vertices.

We see how the differential state calculation procedure described above can be affected by stealthy attacks. For a subset~\({\cal F}\) of edges in~\(G,\) we define a~\({\cal F}-\)\emph{attack vector} or simply an attack vector to be a~\(t \times 1\) vector~\(\mathbf{a} = [a_1,\ldots,a_t]^{T}\) satisfying~\(a_l \neq 0\) if and only if the index~\(l\) corresponds to an edge in~\({\cal F}\) or is a vertex adjacent to an edge in~\({\cal F}.\)

We say that~\({\cal F}\) is \emph{stealthily attackable} if there exists an attack vector~\(\mathbf{a}\) of the form~\(\mathbf{a} = \mathbf{H} \cdot \mathbf{s}\) for some~\(n \times 1\) vector~\(\mathbf{s} = [s_1,\ldots,s_n]^{T}.\) For convenience, we denote~\(\mathbf{a} = \mathbf{H} \cdot \mathbf{s}\) to be a \emph{stealthy} attack vector and denote~\(\mathbf{s}\) to be the~\({\cal F}-\)\emph{stealth vector} or simply the stealth vector corresponding to the attack vector~\(\mathbf{a}.\) In other words, we say that an attack vector~\(\mathbf{a}\) is stealthy if~\(\mathbf{a}\) is also a flow vector.


By definition the attack vector~\(\mathbf{a}\) affects only edges of~\({\cal F}\) or vertices adjacent to edges of~\({\cal F}.\) From the flow equation~(\ref{line_z}), we therefore have that~\(s_i-s_j \neq 0\) if and only if the edge~\((i,j)\) of~\(G\) with endvertices~\(i\) and~\(j\) belongs to~\({\cal F}.\) Given the corrupted flow vector~\(\mathbf{H} \cdot \mathbf{x} + \mathbf{H} \cdot \mathbf{s},\) the differential state calculation procedure described before obtains the difference between the values of the states at vertices~\(i\) and~\(j\) to be
\begin{equation}\label{ph_dif}
\left\{
\begin{array}{cc}
(x_i-x_j) + (s_i-s_j) & \text{ if edge }(i,j) \in {\cal F}\\
(x_i-x_j) & \text{ otherwise}.
\end{array}
\right.
\end{equation}


From~(\ref{ph_dif}), we have that different edges in~\({\cal F}\) are affected differently due to stealthy attacks and so we define the \emph{edge variation factor}~\(\theta({\cal F})\) as
\begin{equation}\label{disp_def}
\theta({\cal F}) := \inf_{\mathbf{s}} \frac{\max_{(i,j) \in {\cal F}} |s_i-s_j|}{\min_{(i,j) \in {\cal F}} |s_i-s_j|},
\end{equation}
where the infimum is taken over all~\({\cal F}-\)stealth vectors~\(\mathbf{s}.\) The edge variation factor measures the variation in the corruption of the state vector~\(\mathbf{x}\) due to stealthy attacks. 

We have the following result regarding the variation factor of stealthy attacks on~\({\cal F}.\) 
\begin{theorem}\label{cap_thm} The set of edges~\({\cal F}\) is stealthily attackable if and only if for every cycle~\(C \in G,\) either~\(C \cap {\cal F} = \emptyset\) or~\(\# C \cap {\cal F} \geq 2.\) Moreover, if~\({\cal F}\) is stealthily attackable, then
\begin{equation}\label{disp_asym}
1 \leq \theta({\cal F}) \leq k({\cal F})-1,
\end{equation}
where~\(k({\cal F})\) is the number of components in the graph~\(G\setminus {\cal F}\) obtained after removing the edges in~\({\cal F}.\)
\end{theorem}
A single edge is stealthily attackable if and only if it is a bridge; i.e., removal of the edge results in disconnection of the graph~\(G\) into two distinct components.

In general, stealthily attacking multiple ``non-critical" edges whose disconnection results in few components, creates a more ``uniform" corruption across the attacked edges. Therefore from the designer perspective, it would be beneficial to have extra protection in these non-critical edges. 







\subsection*{Proof of Theorem~\ref{cap_thm}}
Suppose~\({\cal F}\) is stealthily attackable and there exists a cycle~\(C\) such that\\\(\# C \cap {\cal F} = 1;\) i.e.,
there exists exactly one edge~\(e=(u,v) \in C\) present in~\({\cal F}.\) We arrive at a contradiction as follows.
Let~\(\mathbf{a} = \mathbf{H} \cdot \mathbf{s}\) be the stealthy attack vector on~\({\cal F}\) and let~\(\mathbf{s} = [s_1,\ldots,s_n]^{T}\) be
the corresponding stealth vector. As argued prior to~(\ref{ph_dif}), we have that~\(s_i-s_j \neq 0\) if and only if~\((i,j) \in {\cal F}.\)
Thus~\(s_u \neq s_v.\)

Denoting the cycle~\(C = (u,u_1,u_2,\ldots,u_t,v,u),\) we then have that the edge~\((u,u_1) \notin {\cal F}\)
and so~\(s_u = s_{u_1}.\) Similarly~\((u_1,u_2) \notin {\cal F}\) and so~\(s_{u_1} = s_{u_2}.\) Continuing this way,
we get~\(s_u= s_{u_1} = s_{u_2} =\ldots = s_{u_t}.\) Finally, the edge~\((u_t,v)\) is also not in~\({\cal F}\)
and so~\(s_u = s_{u_t} = s_v,\) a contradiction.

Conversely, suppose every cycle~\(C\) in~\(G\) contains either zero or at least two edges of~\({\cal F}.\) We now use minor reduction to obtain the desired stealth vector, by extracting appropriate components of~\(G\) that allow the vector construction. The details are as follows. First we remove the edges in~\({\cal F}\) to get~\(k = k({\cal F})\) connected components~\({\cal C}_1,\ldots,{\cal C}_k\) of~\(G.\)
Every edge~\(e = (w,y) \in G\) satisfies the following property:
\begin{eqnarray}\label{dis_comp}
&&\text{ The edge~\(e \in {\cal F}\) if and only if the vertices~\(w \in {\cal C}_{w_1}\) and~\(y \in {\cal C}_{y_1}\) } \nonumber\\
&&\text{ belong to distinct components~\({\cal C}_{w_1} \neq {\cal C}_{y_1}.\)}
\end{eqnarray}
Indeed, by construction, we have that if the vertices~\(w\) and~\(y\) belong to distinct components~\({\cal C}_{w_1}\) and~\({\cal C}_{y_1},\) then the edge~\((w,y)\) necessarily belongs to~\({\cal F}.\) Conversely if~\((w,y) \in {\cal F}\) and both~\(w\) and~\(y\) were to belong to the same component~\({\cal C}_j,\) then~\(w\) and~\(y\) would be connected by a path~\(P_{wy}\) in~\({\cal C}_j.\) This in turn would imply that~\(e \cup P_{wy}\) is a cycle in~\(G\) containing only one edge~\(e\) of~\({\cal F},\) a contradiction.

We now construct the stealth vector~\(\mathbf{s} = [s_1,\ldots,s_n]\) as follows. Let~\(\lambda > 0\) be a real number to be determined later. For a vertex~\(v \in G\) we set
\begin{equation}\label{sv_lam}
s_v = s_v(\lambda) = \lambda^{i-1}, \text{ if } v \in {\cal C}_i,\;\;\;\;1 \leq i \leq k
\end{equation}
so that~\(\mathbf{s}(\lambda) = [s_1(\lambda),\ldots,s_n(\lambda)].\) Letting~\(\mathbf{a} = \mathbf{H} \cdot \mathbf{s}\) we first argue that~\(a_l = 0\) if~\(l\) is neither a vertex adjacent to some edge of~\({\cal F}\) nor the index of some edge in~\({\cal F}.\) Indeed if~\(l\) is the index of an edge~\((u,v) \notin {\cal F}\) then both~\(u\) and~\(v\) belong to the same component~\({\cal C}_{i_0}\) for some~\(1 \leq i_0 \leq k\) (property~(\ref{dis_comp})). This implies that~\(s_u = s_v = \lambda^{i_0-1}\) and so from~(\ref{line_z}), we get that~\[|a_l| = B_{u,v}|s_u-s_v| = 0.\] Next if~\(l\) is a vertex and the vertex~\(l\) is not the endvertex of any edge in~\({\cal F},\) then by property~(\ref{dis_comp}) all neighbours of~\(l\) are present in the same component~\({\cal C}_{j_0}\) for some~\(1 \leq j_0 \leq k.\) This implies that~\(s_w = s_l = \lambda^{j_0-1}\) for every neighbour~\(w\) of~\(l\) and from~(\ref{bus_z}), we again get that~\(a_l = 0.\)

We now see that~\(a_l \neq 0\) if either~\(l\) is the index of some edge~\((w,y) \in {\cal F}\) or is a vertex adjacent to some edge of~\({\cal F}.\) First suppose that~\(l\) is the index of~\((w,y)  \in {\cal F}.\) From property~(\ref{dis_comp}) above, we have that~\(w\) and~\(y\) belong to distinct components~\({\cal C}_{w_1} \neq {\cal C}_{y_1}\) and so by construction~\[s_w  = \lambda^{w_1-1} \neq \lambda^{y_1-1} = s_y.\] From~(\ref{line_z}), this implies that~\(|a_l| = B_{w,y}|s_w-s_y| \neq 0.\)

Finally, choosing~\(\lambda\) appropriately, we now show that~\(a_l \neq 0\) if~\(l\) is a vertex adjacent to some edge in~\({\cal F}.\) Suppose the vertex~\(l\) belongs to the component~\({\cal C}_{j_1}\) for some~\(1 \leq j_1 \leq k\) so that the corresponding entry~\(s_l = \lambda^{j_1-1}.\) Since~\(l\) is the endvertex of some edge in~\({\cal F},\) we have from property~(\ref{dis_comp}) that~\(l\) has at least one neighbour outside~\({\cal C}_{j_1}.\) Let~\({\cal C}_{j_1},\ldots,{\cal C}_{j_w}\) be the set of all components containing either~\(l\) or a neighbour of~\(l.\) Using the flow equation~(\ref{bus_z}) we then get that
\begin{equation}\label{zu_eq}
a_l = a_l(\lambda) = \beta_{1,l} \cdot \lambda^{j_1-1} - \sum_{r=2}^{w} \beta_{r,l} \cdot \lambda^{j_{r}-1}
\end{equation}
where~
\[\beta_{1,l} := \sum_{q \sim l} B_{l,q} - \sum_{q \sim l, q \in {\cal C}_{j_1}} B_{l,q} > 0\]
and for~\(2 \leq r \leq w\) the term
\[\beta_{r,l}  := \sum_{q \sim l, q \in {\cal C}_{j_r}} B_{l,q} > 0\]
is the sum of the gains of edges adjacent to the vertex~\(l\) and present in the component~\({\cal C}_{j_r}\) so that
\begin{equation}\label{sum_alph}
\beta_{1,l}  - \sum_{r=2}^{w} \beta_{r,l} = 0.
\end{equation}


Let~\(\Lambda_l\) be the finite set of the roots of the equation~\(a_l(x) = 0\) so that~\(1 \in \Lambda_l\) by~(\ref{sum_alph}) and set~\(\Lambda_{tot} = \bigcup_v \Lambda_v,\)  where the union is taken over all vertices adjacent to an edge in~\({\cal F}.\) Choosing~\(\lambda \notin \Lambda_{tot}\) we have that~\(a_l(\lambda) \neq 0\) and since~\(l\) is arbitrary this implies that~\(\mathbf{a}(\lambda) = \mathbf{H} \cdot \mathbf{s}(\lambda)\) is a stealthy attack vector.

From~(\ref{sv_lam}) and property~(\ref{dis_comp}), we also get for any edge~\((i,j) \in {\cal F}\) that the corresponding difference~\(|s_i(\lambda)-s_j(\lambda)|\) is of the form~\(|\lambda^{i_1}-\lambda^{j_1}|\) for some distinct integers~\(i_1,j_1 \geq 0.\) Therefore if~\(\{\lambda_q\}_{q \geq 1}, \lambda_q  \notin \Lambda_{tot}, \lambda_q < 1\) is a sequence converging to one as~\(q \rightarrow \infty\) then using
\[\max_{(i,j) \in {\cal F}} |s_i(\lambda_q)-s_j(\lambda_q)| = 1-\lambda_q^{k-1} \text{ and } \min_{(i,j) \in {\cal F}} |s_i(\lambda_q)-s_j(\lambda_q)| = \lambda_q^{k-2}-\lambda_q^{k-1}\] we get that
\begin{equation}\label{disp_est}
\frac{\max_{(i,j) \in {\cal F}} |s_i(\lambda_q)-s_j(\lambda_q)|}{\min_{(i,j) \in {\cal F}} |s_i(\lambda_q)-s_j(\lambda_q)|} = \frac{1-\lambda_q^{k-1}}{\lambda_q^{k-2}(1-\lambda_q)} \longrightarrow k-1
\end{equation}
as~\(q \rightarrow \infty.\) This implies that~\(\theta({\cal F}) \leq k({\cal F})-1\)
and completes the proof of Theorem~\ref{cap_thm}.~\(\qed\)



\setcounter{equation}{0}
\renewcommand\theequation{\thesection.\arabic{equation}}
\section{Extension of results in Theorem~\ref{cap_thm}} \label{prac_sec}
In this section, we extend the result of Theorem~\ref{cap_thm} in two directions.
In the first subsection, we use the structure of the graph~\(G \setminus {\cal F}\) to improve the bound for the edge variation factor and in the second subsection, we study stealthy attacks with partial information.

\subsection*{Improved Bound for the Variation Factor}
We now describe a slight modification of the proof of Theorem~\ref{cap_thm} to obtain stealth vectors with lesser variation. Consider the component graph~\(G_{comp}=  (V_{comp}, E_{comp})\) constructed as follows: We represent the component~\({\cal C}_i\) by a node~\(n_i \in V_{comp}\) and connect nodes~\(n_i\) and~\(n_j\) if there exists an edge~\(e \in {\cal F}\) with one endvertex in~\({\cal C}_i\) and other endvertex in~\({\cal C}_j.\)
A proper colouring of~\(G_{comp}\) using~\(a \geq 1\) colours is a map~\(g : V_{comp} \rightarrow \{1,2,\ldots,a\}\) such that~\(g(v) \neq g(u)\) if vertices~\(u\) and~\(v\) are adjacent in~\(G_{comp}.\) The chromatic number~\(\chi_0\) of~\(G_{comp}\) is the smallest integer~\(z\) such that~\(G_{comp}\) admits a proper colouring using~\(z\) colours~\cite{boll}.

Letting~\(g_0 : V_{comp} \rightarrow \{1,2,\ldots,\chi_0\}\) be a proper colouring of~\(G_{comp}\) using~\(\chi_0\) colours, we define the stealth vector~\(\mathbf{y} = [y_1,\ldots,y_n]^{T}\) as follows. If~\(g_0(n_i) = j,\) we assign~\(y_v = y_v(\lambda) = \lambda^{j-1}\) for all vertices~\(v \in {\cal C}_i.\) Finally, letting~\(\mathbf{y} = [y_1,\ldots,y_n]^{T},\) we set~\(\mathbf{c} = \mathbf{H} \cdot \mathbf{y}.\)

To see that~\(\mathbf{c}\) is a stealthy attack vector, we consider any edge~\((u,v) \in G\) with index~\(l.\) If~\(u\) and~\(v\) belong to the same component in~\(\{{\cal C}_i\},\) then~\(y_u-y_v = 0\) and so we have from~(\ref{line_z}) that the corresponding entry~\(c_l = B_{u,v}(y_u-y_v) =0.\) Similarly if a vertex~\(u\) is not adjacent to any edge in~\({\cal F}\) then using property~(\ref{dis_comp}) as before, we get that all neighbours of~\(u\) belong to the same component in~\(\{{\cal C}_i\}\) as~\(u.\) As before, we use~(\ref{bus_z}) to get that the corresponding entry~\(c_u = 0.\)

If an edge~\((u,v) \in {\cal F}\) then by property~(\ref{dis_comp}), the vertices~\(u \in {\cal C}_{u_1}\) and~\(v \in {\cal C}_{v_1}\) belong to distinct components~\({\cal C}_{u_1} \neq {\cal C}_{v_1}.\) In the graph~\(G_{comp}\) constructed above, the nodes~\(n_{u_1}\) and~\(n_{v_1}\) are joined by an edge and so have different colours~\(g(n_{u_1}) \neq g(n_{v_1}).\) This in turn implies that~\(y_u \neq y_v\) and so~\(c_l = B_{u,v}(y_u-y_v) \neq 0,\) by~(\ref{line_z}).

Finally, for a vertex adjacent to an edge in~\({\cal F},\) we argue as in~(\ref{zu_eq}) to get the corresponding polynomial expression for~\(c_l(\lambda).\) The expressions here have different degrees than in~(\ref{zu_eq}) and therefore different set of roots. Again we choose an appropriate sequence~\(\{\lambda_q\}\) converging to one so that
\[\frac{1-\lambda_q^{\chi_0-1}}{\lambda_q^{\chi_0-2}-\lambda_q^{\chi_0-1}} \longrightarrow \chi_0-1\] as~\( q \rightarrow \infty.\) This implies that~\(\theta({\cal F}) \leq \chi_0-1.\)

\subsection*{Partial Knowledge}
In this subsection, we see how stealth attacks could possibly be carried out with coarse information regarding the gain matrix~\(\mathbf{H}.\) This situation arises for example, when measurement noise results in imperfect estimates of the gain values.

Suppose the attacker does not exactly know the true gains~\(\{B_{i,j}\}\) but knows that
\begin{equation}\label{eps_bound}
\epsilon_1 \leq B_{i,j} \leq \epsilon_2
\end{equation}
for some positive finite constants~\(\epsilon_1,\epsilon_2.\)

We construct the stealth vector~\(\mathbf{s}\) as follows. For a vertex~\(v \in G,\) we choose~\(s_v\) as in~(\ref{sv_lam}) and get from the proof of Theorem~\ref{cap_thm} that~\(a_l = 0\) if~\(l\) is neither the index of an edge in~\({\cal F}\) nor is a vertex adjacent to any edges in~\({\cal F}.\) Moreover~\(a_l \neq 0\) if~\(l\) is the index of an edge of~\({\cal F}.\) The knowledge of edge gains is required only to determine an appropriate value of~\(\lambda\) so that~\(a_l \neq 0\) for a vertex~\(l\) adjacent to some edge in~\({\cal F}.\) Indeed from~(\ref{zu_eq}), we have that
\begin{equation}\label{zu_eq2}
a_l = a_l(\lambda,\mathbf{H}) = \beta_{1,l}\cdot \lambda^{j_1-1} - \sum_{r=2}^{w} \beta_{r,l} \cdot \lambda^{j_{r}-1}
\end{equation}
where the positive ~\(\{\beta_{r,l}\}\) are as in~(\ref{zu_eq}).

We now see that if~\(\lambda > 0\) is chosen large, then~\(|a_l| \geq D\) for some constant\\\(D = D(\lambda,\epsilon_1,\epsilon_2)\) not depending on the choice of~\(\mathbf{H}.\) Indeed, assuming\\\(j_1>j_2\ldots > j_w\) in~(\ref{zu_eq2}) and using~(\ref{eps_bound}) we have
\begin{eqnarray}
|a_l| &=&  \beta_{1,l}\lambda^{j_1-1} \left| 1- \sum_{r=2}^{w} \frac{\beta_{r,l}}{\beta_{1,l} \cdot \lambda^{j_1-j_r}} \right| \nonumber\\
&\geq& \epsilon_1 \lambda^{j_1-1}\left(1- \sum_{r=2}^{w} \frac{\epsilon_2}{\epsilon_1 \cdot \lambda^{j_1-j_r}}\right) \nonumber
\end{eqnarray}
so that
\begin{equation}\label{al_est}
|a_l| \geq \epsilon_1 \lambda^{j_1-1}\left(1- \frac{w \cdot \epsilon_2}{\epsilon_1 \cdot \lambda^{j_1-j_r}}\right) \geq  \epsilon_1 \lambda^{j_1-1}\left(1- \frac{k \cdot \epsilon_2}{\epsilon_1 \cdot \lambda^{j_1-j_r}}\right)
\end{equation}
where~\(k = k({\cal F})\) is the number of components of~\(G\setminus {\cal F}.\) The final relation in~(\ref{al_est}) is true because, the number of distinct powers of~\(\lambda\) in the stealth vector~\(\mathbf{s}\) is~\(k.\)  Thus for all~\(\lambda \geq \lambda_0(l,k,\epsilon_1,\epsilon_2) > 1\) large we have from~(\ref{al_est})
that~\[|a_l| \geq \frac{\epsilon_1 \lambda^{j_1-1}}{2} \geq \frac{\epsilon_1}{2}.\] Choosing~\(\lambda > \lambda_1 = \lambda_1(k,\epsilon_1,\epsilon_2)\) large enough, we get that~\(|a_l| \geq \frac{\epsilon_1}{2}\) for \emph{all} vertices~\(l\) adjacent to some edge in~\({\cal F}.\) This obtains our desired stealth vector~\(\mathbf{s} = \mathbf{s}(\lambda).\)

As a concluding remark, we state here that though the attack as described is possible, the edge variation~\(\theta({\cal F})\) could be quite large and therefore be possibly detected.~\(\qed\)

\setcounter{equation}{0}
\renewcommand\theequation{\thesection.\arabic{equation}}
\section{Conclusion and Future Work} \label{sec_conc}
In this paper, we have studied two applications of graph minor reduction in boxicity and stealthy attacks on flowgraphs. We have used minor reduction recursively to estimate the boxicity of an arbitrary graph. Similarly, we used a component reduction technique to determine necessary and sufficient conditions that allows stealthy attacks in a flowgraph.

In the above, we have considered deterministic graphs. In the future, we plan to incorporate and study analogous problems in random graphs.

\begin{acknowledgments}
  I thank Professors V. Raman, C. R. Subramanian and the referee for crucial comments that led to an improvement of the paper. I also thank IMSc for my fellowships.
\end{acknowledgments}

\bibliographystyle{plain}


\end{document}